\documentclass{amsart}

\usepackage{amssymb,latexsym}
\usepackage{amsmath}
\usepackage{graphicx}
\usepackage{textcomp}

\usepackage{amsthm,amssymb,enumerate,graphicx, tikz}
\usepackage{amscd}
\usepackage{setspace}
\usepackage{comment}
\usepackage{hyperref}
\usepackage{cleveref}

\usepackage{xcolor}

\title[rainbow versions of the Caccetta-H\"aggkvist conjecture]{Non-uniform degrees and rainbow versions of the Caccetta-H\"aggkvist conjecture}

\author{Ron Aharoni} \thanks{We acknowledge the financial support from the Ministry of Educational and Science of the Russian Federation in the framework of MegaGrant no. 075-15-2019-1926 when the first author worked on \refS{sec:triangles} of the paper.\\
\indent Ron Aharoni:
Department of Mathematics, Technion, Israel and MIPT. \url{ra@technion.ac.il}. The research of R.~Aharoni was supported in part by the Israel Science Foundation (ISF) grant no. 2023464 and the Discount Bank Chair at the Technion. This paper is part of a project that has received funding from the European Union's Horizon 2020 research and innovation programme under the Marie Skldowska-Curie grant agreement no.\ 823748.}
\author{Eli Berger} \thanks{Eli Berger: Department of Mathematics, University of Haifa,  Israel. \url{berger@math.haifa.ac.il}.}
\author{Maria Chudnovsky} \thanks{Maria Chudnovsky: Department of Mathematics, Princeton University, USA, \url{mchudnov@math.princeton.edu}.}
\author{He Guo}\thanks{He Guo: Faculty of Mathematics, Technion, Israel. \url{hguo@campus.technion.ac.il}.}
\author{Shira Zerbib} \thanks{Shira Zerbib: Department of Mathematics, Iowa State University, USA.  \url{zerbib@iastate.edu}. The research of S.~Zerbib was supported by NSF grant DMS-1953929.\\
\indent The authors were supported by US-Israel Binational Science Foundation (BSF) grant no.\ 2016077.}
\newtheorem{theorem}{Theorem}[section]

\newtheorem*{theorem*}{Theorem}

\newtheorem{question}[theorem]{Question}
\newtheorem{problem}[theorem]{Problem}

\newtheorem{corollary}[theorem]{Corollary}

\newtheorem{claim}{Claim}[theorem]
\newtheorem{conjecture}[theorem]{Conjecture}
\theoremstyle{definition}
\newtheorem{definition}[theorem]{Definition}

\theoremstyle{remark}
\newtheorem{remark}[theorem]{Remark}
\newtheorem{example}[theorem]{Example}

\newcommand{\refT}[1]{Theorem~\ref{#1}}
\newcommand{\refC}[1]{Corollary~\ref{#1}}

\newcommand{\refS}[1]{Section~\ref{#1}}

\newcommand{\refCl}[1]{Claim~\ref{#1}}
\newcommand{\refCon}[1]{Conjecture~\ref{#1}}

\newcommand{\cPr}{\mathbb{P}}
\newcommand{\E}{\mathbb{E}}
\newcommand{\Bin}{\operatorname{Bin}}
\newcommand{\cS}{\mathcal{S}}
\newcommand{\cA}{\mathcal{A}}
\newcommand{\cB}{\mathcal{B}}
\newcommand{\cC}{\mathcal{C}}

\newcommand{\cf}{\mathcal{F}}

\newcommand{\dist}{\text{dist}}
\begin{document}

\begin{abstract}
%The  Caccetta-Haggkvist conjecture states that in any  directed graph on $n$ vertices with all  out-degrees  at least $r$ there is a directed cycle of size at most $\lceil n/r \rceil$. The conjecture is proved for $r\le 5$ and ... The first author conjectured a rainbow generalization of this conjecture: if $G$ is a non-directed graph 

The Caccetta-H\"aggkvist  conjecture (denoted below  CHC) states that the directed girth (the smallest length of a directed cycle) $dgirth(D)$ of a directed graph $D$ on $n$ vertices is at most $\lceil \frac{n}{\delta^+(D)}\rceil$, where $\delta^+(D)$ is the minimum out-degree of~$D$.  
We consider a version involving all out-degrees, not merely the minimum one, and prove that if $D$ does not contain a sink, then $dgirth(D) \le 2 \sum_{v\in V(D)} \frac{1}{deg^+(v)+1}$. 
In the spirit of a generalization of the CHC to rainbow cycles  in \cite{ADH2019}, this  suggests the conjecture that given non-empty sets $F_1, \ldots,F_n$ of edges of $K_n$, there exists a rainbow cycle of length at most $2\sum_{1\le i \le n}\frac{1}{|F_i|+1}$. We prove a bit stronger result when $1\le |F_i|\le 2$, thereby 
strengthening a result of DeVos et. al \cite{DDFGGHMM2021}.  
We prove a  logarithmic bound on the rainbow girth in the case that the  sets $F_i$ are triangles. 

%We prove a special case of a rainbow version of the CH conjecture (suggested in  \cite{ADH2019}: given $n$ non-empty sets $F_i$ of edges in an undirected  graph on $n$ vertices,  there exists a rainbow cycle of length no larger than $\lceil \sum_{v\in V(G)} 2\frac{1}{|F_i|+1}\rceil$. We prove this in the case that $|F_i| \le 2$ for all $i\le n$. 

%%,  we conjecture that every digraph $G$ has a directed cycle of size at most $\lceil\sum_{v\in V(G)} \frac{1}{deg^+(v)}\rceil$. 
%We prove this when $deg^+(v)\le 2$ for all $v$.
%We also address a rainbow generalization of the conjecture: if $F_1,\dots,F_n \subseteq E(K_n)$, then they have a rainbow cycle  of size  $\lceil \sum_{i=1}^n \frac{1}{|F_i|}\rceil$ or less. (``Rainbow'' means choosing a distinct representative for some of the sets. The CH conjecture is the case in which  the $F_i$s are the out-stars in the directed graph.) We prove this conjecture when $|F_i| \le 2$ for all $i$. 

\end{abstract}

\maketitle

\section{Introduction}\label{sec:introduction}
The {\em directed girth} $dgirth(D)$ of a directed graph (digraph)~$D$ is the smallest length of a directed cycle in~$D$ ($\infty$ if there is no directed cycle). A famous conjecture of Caccetta and H\"aggkvist~\cite{CaccettaHaggkvist} is that 
 \[dgirth(D)\le \Big\lceil \frac{n}{\delta^+(D)}\Big\rceil,\]
 where $n=|V(D)|$ and 
 $\delta^+(D)$ is the minimum out-degree over all vertices of $D$. We use the acronym CHC for it. See~\cite{sullivan} for a survey of known results on this conjecture up to the year 2006. 
 
The CHC is known to be true asymptotically: in~\cite{shen3}
it was proved that 
 \begin{equation}\label{asymptotic}
 dgirth(D) \le \Big\lceil \frac{n}{\delta^+(D)}\Big\rceil +73. 
 \end{equation}
 
 Much of the research on the conjecture has addressed the case $dgirth(D)=3$.
 %$\delta^+(D)=\frac{n}{3}$. 
 The best result so far is due to Hladk{\'y}, Kr{\'a}l', and Norin~\cite{hkn}.
 
\begin{theorem}\label{3465}
 Every $n$-vertex digraph with minimum out-degree at least 0.3465n contains a directed triangle.
\end{theorem}

 A natural question is finding upper bounds on $dgirth(D)$ in terms of all out-degrees of the vertices of~$D$, rather than merely the minimum out-degree. 
Let 
 \[\psi(D):=\sum_{v \in V(D)}\frac{1}{\deg^+(v)}. \]

Seymour asked (see \cite{hompe}) whether CHC could be generalized to 
 \begin{equation}\label{seymoureq}
 dgirth(D) \le \lceil \psi(D) \rceil.
 \end{equation}
This was answered in the negative  by Hompe~\cite{hompe}. 
Here we  prove ``half'' of this result, namely:

\begin{theorem}
For any digraph~$D$, we have
\begin{equation}\label{eq:twiceseymoureq}
 dgirth(D) \le 2 \psi(D).
 \end{equation}

\end{theorem}

In fact, we use a slightly different function. Let

\[ \varphi(D):= 
\sum_{v\in V(D)} \frac{1}{deg^+(v)+1}. \]

\begin{theorem}\label{main}
If all out-degrees  in $D$ are positive, then 
$ dgirth(D) \le 2\varphi(D)$. 
\end{theorem}

%\textcolor{red}{Also, what is the connection of  Question 1.2 and Thm 1.6 below? there it is written that the result is sharp, with a constant 2, and here the constant is less than 2. Please clear this point}

This is proved in~\refS{sec:2}.

In~\refS{sec:triangles} and~\ref{sec:max2} we discuss a rainbow, undirected generalization of the CHC, suggested in \cite{ADH2019}.

\begin{definition} 
Let $\cf=(F_1, \ldots, F_m)$ be a family of subsets of $E(K_n)$.
A \emph{rainbow cycle} for~$\cf$ is a cycle whose edges are chosen each from a different $F_i$. The {\em rainbow girth} $rgirth(\cf)$ of $\cf$ is the smallest length of a rainbow cycle. 
\end{definition}
Note that an edge belonging to two different sets ~$F_i$  yields a rainbow digon (that is, a rainbow cycle of length 2). Thus for our purposes
we can assume disjointness of the sets $F_i$. The generalized CHC is:

\begin{conjecture}\label{con:generalCHC}
For $\cf=(F_1, \ldots, F_n)$ a family of subsets of $E(K_n)$, we have $rgirth(\cf) \le \lceil \frac{n}{\min_{1\le i\le n} |F_i|}\rceil$.
\end{conjecture}

As explained in \refS{sec:triangles}, the CHC is the case in which the sets $F_i$ are stars, with distinct apexes. 

\begin{remark}
An advantage of the rainbow version is that it detaches the link between the number of sets and the number $n$ of vertices. The question makes sense for any number of sets. Here are two results on the case $rgirth(\cf)=3$:

\begin{theorem}\cite{GoorevitchHolzman}
$n^2/8+o(n)$ triangles on $n$ vertices have a rainbow triangle.
\end{theorem}

\begin{theorem}\cite{ADH2019}
$\frac{9}{8}n$ (or more) sets of edges in $K_n$, each of size $\frac{n}{3}$ or more,  have a rainbow triangle.
\end{theorem}

In \cite{HQS2022} a slight improvement was proved, $\frac{9}{8}n$ being replaced by $1.1077n$.
\end{remark}

In~\cite{HS2022} it was shown that the order of magnitude in the conjecture is correct: 

\begin{theorem}\label{constant}
There exists a constant $0< C\le 10^{11}$ such that for any~$n$ and any family $\cf=(F_1, \ldots, F_n)$ of subsets of~$E(K_n)$,  we have $rgirth(\cf) \le C\cdot \frac{n}{\min_{1\le i\le n} |F_i|}$.
\end{theorem}

A natural challenge is to improve the bound on~$C$. 
%\\ \ \\

In \cite{ADH2019} a triangles version was proved: 
\begin{theorem}
$n$ sets of edges in $K_n$, each of size $0.4n$ or more,  have a rainbow triangle.
\end{theorem}
Compare with the coefficient $0.3465$ appearing in Theorem \ref{3465}. In~\cite{HQS2022},~$0.4n$ was replaced by $0.3988n$.

%\\ \ 

In \cite{DDFGGHMM2021} the following was proved:
\begin{theorem}\label{devos}
\refCon{con:generalCHC} is true when $|F_i|= 2$ for all $1\le i \le n$.
\end{theorem}

The rainbow analogue of \refT{main} is:

\begin{conjecture}\label{rainbowmain}
$rgirth(\cf) \le 2\sum_{1\le i \le n}\frac{1}{|F_i|+1}$ for any family $\cf=(F_1,\dots, F_n)$ of subsets of~$E(K_n)$.
\end{conjecture}
If true, this would enable taking $C=2$ in \refT{constant}. Here we prove:
\begin{theorem}
\refCon{rainbowmain} is true when $1\le |F_i|\le 2$ for all $1\le i \le n$.
\end{theorem}
We shall prove this (in ~\refS{sec:max2}) via a result  generalizing~\refT{devos}. 
Let 
\[\psi(\cf) := \sum_{1\le i \le n}\frac{1}{|F_i|}.\] 
We show: 
\begin{theorem}\label{rainbowpsi}
If $|F_i|\le 2$ for all $1\le i \le n$, then 
$rgirth(\cf) \le \lceil \psi(\cf) \rceil $.
\end{theorem}
\refS{sec:triangles} deals with a special case of \refCon{rainbowmain}, in which all sets~$F_i$ are triangles. This  case is of particular interest, for the following reason. We know (from the original CHC) that $\min |F_i|\cdot rgirth(\cf)$ may be close to~$n$, and that this can be exactly~$n$ when the sets~$F_i$ are stars. In \cite{AharoniGuo} it was proved that if each $F_i$ is a matching of size $2$ then $rgirth(\cf) =O(\log n)$. Note that a set of
%at least two 
%{\color{red} a single edge is also a star}
graph edges  not containing two disjoint edges is a star or a triangle. So, the remaining case, in terms of some uniform assumption on the sets $F_i$, is that of triangles. We  show that this case is close to the case of matchings of size $2$:

\begin{theorem}
    For any constant $\alpha>1/2$ there exists a constant~$C$ such that for any~$n$ and any family $\cf=(F_1,...,F_{\lceil \alpha n\rceil})$ of subsets of~$E(K_n)$ where each $F_i$ is a triangle, there is a rainbow cycle of length at most~$C\log n$.
\end{theorem}

We also prove, via a random construction, that this result is best possible, in the sense that there are families of $n$ triangles, in which the rainbow girth is~$\Omega(\log n)$. (For a stronger version, see~\refT{thm:nosmallcycle}.)
\begin{theorem}\label{Omegalogn}
There exists a positive constant~$c$ such that for any~$n$, there exists a family~$\cf_n$ of~$n$ triangles on~$n$ vertices satisfying $rgirth (\cf_n) \ge c \log n$.
\end{theorem}

\section{Non-uniform out-degrees}\label{sec:2}
As mentioned in~\refS{sec:introduction}, Seymour asked (see \cite{hompe}) whether the directed girth of a digraph can be bounded from above by an expression involving all out-degrees. A natural such expression is
\[\psi(D)=\sum_{v \in V(D)}\frac{1}{\deg^+(v)}. \]

Hompe~\cite{hompe} showed that $\psi(D)$ is not always an upper bound on the directed girth.
His counterexample is obtained from a directed cycle by replacing each vertex of the cycle by  a transitive tournament $T_k$ with $k$ vertices, for some $k$. 
If the cycle is of length $\ell$ then the resulting graph~$D$ satisfies $dgirth(D)=\ell$ and $\delta^+(D)=k$. Furthermore, $\psi(D)=dgirth(D)\cdot \sum_{i=\delta^+(D)}^{2\delta^+(D)-1}\frac{1}{i}$, and  $\varphi(D)=dgirth(D)\cdot \sum_{i=\delta^+(D)}^{2\delta^+(D)-1}\frac{1}{i+1}$, and 
$\lim_{|V(D)|\to \infty} \frac{dgirth(D)}{\psi(D)}= \lim_{|V(D)|\to \infty} \frac{dgirth(D)}{\varphi(D)} = \log_2 e$. 

Possibly, this example is best: 
\begin{question}\label{new}
Is it true that for any digraph~$D$,
$dgirth(D) \le \lceil \log_2 e \cdot \psi(D) \rceil$? 
\end{question}

 A {\em sink} in a digraph is a vertex with out-degree $0$. If $D$ contains a sink, then $\psi(D)=\infty$, and thus $dgirth(D) \le \psi(D)$. Thus the interesting case for us is that in which no sink exists. In this case we can prove twice the bound suggested by Seymour, in fact a bit better. Recall that
\[\varphi(D)=\sum_{v \in V(D)}\frac{1}{\deg^+(v) +1}.\]

\begin{theorem}\label{thm:dgphi}
If a digraph~$D$ has no sink, then
$dgirth(D) \le 2 \varphi(D)$. Equality holds if and only if $D$ is a Hamilton cycle (in which $dgirth(D) = 2\varphi(D) =|V(D)|$) or a complete digraph (in which case $dgirth(D) = 2\varphi(D) =2$).
\end{theorem}

In~\cite{chvatalszemeredi} the inequality was proved in the case that all out-degrees are equal. 
% The inequality is probably far from sharp in other cases than Hamilton cycles and complete digraphs. 

% For example, in  Hompe's construction when the minimum out-degree is $2$ , $dgirth(D)= \frac{12}{7} \cdot \varphi(D) $. 
%\begin{question}\label{new:delta=2}
%Is it true that  if $\delta^+(D) \ge 2$ then $dgirth(D) \le \lceil \frac{12}{7} \cdot \varphi(D) \rceil$?
%\end{question}

\begin{proof}[Proof of \refT{thm:dgphi}]     
Let us first prove the inequality. 
We call a digraph $K$ not containing a sink {\em $\varphi$-critical}
if for every vertex $v \in V(K)$ either $\varphi(K-v)>\varphi(K)$ or $K-v$ contains a sink.

\begin{claim}\label{claim:phicritical}
A $\varphi$-critical graph is vertex-disjoint union of directed cycles. \end{claim}
 
\begin{proof}[Proof of~\refT{thm:dgphi} based on~\refCl{claim:phicritical}]
We remove vertices one by one from~$D$, while
keeping the graph sink-less and not increasing~$\varphi$, until we reach a $\varphi$-critical graph~$K$ that is vertex-disjoint union of directed cycles. Since~$K$ is
union of cycles, we have $dgirth(K)\le |V(K)|=2\varphi(K)$.
Since~$K$ is a subgraph of~$D$, we have $dgirth(D) \le dgirth(K)$. Since we keep~$\varphi$ not increasing during the removal, we have $\varphi(K)\le \varphi(D)$. Combining these, we have
\[ dgirth(D)\le dgrith(K)\le 2\varphi(K)\le 2\varphi(D), \]
which completes the proof. 
\end{proof} 

To prove~\refCl{claim:phicritical} we observe:

\begin{claim}
In any digraph~$D$, there exists a vertex $v$ for which 
\begin{equation} \label{eq:claim}
\frac{1}{deg^+(v)+1} \ge 
\sum_{u \in N^-(v)}\frac{1}{deg^+(u)}\frac{1}{deg^+(u)+1}\end{equation}
\end{claim}

\begin{proof}
The claim will follow if we show that 
the sums, over all vertices of~$D$,  of the two sides, are the same. On the left-hand side the sum is, by definition, $\varphi(D)$. On the right-hand side, the number of times every vertex~$u$ appears is $deg^+(u)$, and hence we get $\sum_{u \in V(D)} \frac{1}{deg^+(u)+1}$, which is again~$\varphi(D)$.
\end{proof}

\begin{proof}[Proof of~\refCl{claim:phicritical}]
Let~$D$ be a~$\varphi$-critical graph and~$A$ be the set of vertices~$v$ satisfying~\eqref{eq:claim}.
Note that for any $v\in A$, 
\[ \varphi(D)-\varphi(D-v)= \frac{1}{deg^+(v)+1}-
\sum_{u \in N^-(v)}\Big(\frac{1}{deg^+(u)}-\frac{1}{deg^+(u)+1}\Big)\ge 0. \]
As~$D$ is~$\varphi$-critical,  for every~$v\in A$, $D-v$ has a sink, which means there exists a vertex~$w$ such that $N^{+}(w)=v$.
Then the $w$-term in the right-hand side of \eqref{eq:claim} is $\frac{1}{2}$, while the left-hand side is at most $\frac{1}{2}$ as~$D$ is sink-less, and thus $N^{-}(v)=\{w\}$ and $deg^+(v)=1$. Namely, both in-degree and out-degree of~$v$ are $1$. 
It follows that for every $v \in A$ equality holds in \eqref{eq:claim}, and since the sums over all vertices~$v$ of the right-hand sides and the left-hand sides in~\eqref{eq:claim} are equal, 
it implies that~$A=V(D)$. Therefore every vertex of~$D$ has both in-degree and out-degree equal to~$1$, which means~$D$ is vertex-disjoint union of directed cycles. \end{proof}

This concludes the proof of the inequality in \refT{thm:dgphi}.

%\begin{theorem}\label{thm:dg=phi}
  %  If a digraph $D$ has no sink and satisfies $dgirth(D)=2\varphi(D)$, then~$D$ is either a Hamilton directed cycle or a complete directed graph.  
%\end{theorem}

For the second part of the theorem, assume that $dgirth(D)=2\varphi(D)$. Tracking the proof of the inequality,  for $0\le i\le t$ let $D_i=D-\{v_j\mid 1\le j \le i\}$ (so $D_0=D$), where $v_1, \ldots, v_t$ are the removed vertices  from~$D$ (if any), in this order. Then
\[dgirth(D_{i-1}) \le dgirth(D_{i}) \le 2\varphi (D_i) \le 2\varphi(D_{i-1}),    \]
where the second inequality is by the first part of this theorem. 
By the assumption that $dgirth(D)=2\varphi(D)$, equalities hold throughout, namely $dgirth(D_i) = dgirth(D_{i-1})=2\varphi(D_i)$  and $\varphi(D_i) = \varphi(D_{i-1})$. Let $K=D_t$. By the construction and~\refCl{claim:phicritical}, the~$\varphi$-critical graph~$K$ is the vertex-disjoint union of directed cycles, and since $dgirth(K)=2\varphi(K)$ it is a single cycle, namely it is a Hamilton cycle.

If $V(K)=V(D)$, then ~$D$ itself is a Hamilton cycle, proving the desired result. So, we can assume that  $V(K)\subsetneqq V(D)$. 

\begin{claim} 
If $V(K)\subsetneqq V(D)$, then $K$ is a directed 2-cycle, i.e., a directed digon. 
\end{claim} 

To show this, let $p=|N_{D_{t-1}}^+(v_t)\cap V(K)|$ and $q=|N_{D_{t-1}}^-(v_t)\cap V(K)|$. Then the fact that ~$\varphi(D_{t-1})=\varphi(K)$ 
 implies that 
    \[ \frac{1}{p+1}=q\Big(\frac{1}{2}-\frac{1}{3}\Big)=\frac{q}{6}. \]
    Therefore we have $(p,q)=(5,1),(2,2)$, or $(1,3)$.
Since  $dgirth(K)=dgirth(D_{t-1})$  we have $(p,q)=(2,2)$ and~$K$ is a digon, otherwise~$D_{t-1}$ has a shorter directed cycle than~$K$. This proves the claim, and implies that   $D_{t-1}$ is a complete directed graph on three vertices.

This was the first step in the inductive proof of the following claim: 

\begin{claim}
If $V(K)\subsetneqq V(D)$, then $D_i$ is the complete digraph on $2+t-i$ vertices for all $0\le i \le t$.
\end{claim} 

We prove this by induction on $|V(D)|-i$. Assuming that  $D_i$ is complete digraph on $2+t-i$ vertices, let  
$p=|N_{D_{i-1}}^+(v_{i})\cap V(D_{i})|$ and $q=|N_{D_{i-1}}^-(v_{i})\cap V(D_i)|$. Since $\varphi(D_i)=\varphi(D_{i-1})$, we have
\[ \frac{1}{p+1}= q\Big(\frac{1}{|V(D_{i})|}-\frac{1}{|V(D_{i})|+1}\Big)=\frac{q}{|V(D_{i})|(|V(D_{i})|+1)}. \]
Since  $0\le p,q\le |V(D_{i})|$, we have $p=q=|V(D_{i})|$, so  $D_{i-1}$ is a complete digraph on~$|V(D_{i})|+1=2+t-(i-1)$ vertices. This completes the proof of the claim.

Putting $i=0$ proves the statement in the theorem. 
\end{proof}

\section{The rainbow version of CHC for triangles } \label{sec:triangles}

In this section and the next we consider the rainbow, undirected generalization of the CHC.

Here is an explanation why \refCon{con:generalCHC} is a generalization of CHC.
For a directed edge $e=(u,v)$  let $n(e)$
be the undirected pair $\{u,v\}$. Given a digraph~$D$, for every vertex~$u\in V(D)$ let $F_u=\{n(uv) \mid (u,v) \in E(D)\}$  be the star  of edges leaving $u$,  with their direction removed.  Let~$G(D)$ be an undirected graph with vertex set~$V(D)$ and edge set~$\cup_{u\in V(D)}F_u$. Note that sets~$F_u$ are stars with distinct apexes in~$G$. It is easy to verify that a sequence of vertices $v_1v_2\ldots v_k$ forms a rainbow cycle in~$G$ if and only if they form a directed cycle in~$D$.

The CHC holds asymptotically: it is known that $dgirth(D) \le \lceil \frac{n}{\delta^+(D)}\rceil +73$ (see~\cite{shen3}). In the undirected rainbow version the gap between the conjecture and the known bounds is much larger.

In~\cite{AharoniGuo} it was proved that there exists a constant~$C$ for which every set of~$n$ matchings of size $2$ in $K_n$ has a rainbow cycle of length at most $C\log n$. If~$\cf=(F_1,\dots, F_n)$ are~$n$ stars with distinct apexes then directing all edges in $F_i$ away from the apex yields, by~\refT{thm:dgphi}, we have that $rgirth(\cf) \le 2\psi(\cf)$. We cannot prove the same if the apexes are allowed to coincide:

\begin{problem}
Prove (or disprove) $rgirth(\cf) \le 2\psi(\cf)$ for any set of~$n$ stars in~$K_n$.
\end{problem}

Since a set of edges  not containing a matching of size $2$ is either a star  or a triangle,  the remaining case (assuming all sets~$F_i$ are of size at least $2$) is that of triangles. Like in the case of sets of edges containing each a pair of disjoint edges, a better than linear bound can be proved in this case:

\begin{theorem}\label{thm:triangle}
    For any constant $\alpha>1/2$ there exists a constant~$C$ such that for any~$n$ and any family $\cf=(F_1,...,F_{\lceil \alpha n\rceil})$ of subsets of~$E(K_n)$ where each $F_i$ is a triangle, there is a rainbow cycle of length at most~$C\log n$.
\end{theorem}

The proof uses the following result of Bollob\'as and Szemer\'edi~\cite{BS2002}.

\begin{theorem}\label{thm:girthsparsegraph}
For $n\ge 4$ and $k\ge 2$, every $n$-vertex graph with $n+k$ edges has girth at most
\[\frac{2(n+k)}{3 k}(\log k+\log\log k+4). \]
\end{theorem}

\begin{proof}[Proof of~\refT{thm:triangle}]
As noted above, we may assume that the sets $F_i$ are edge-disjoint, or else $rgirth(\cf)=2$.
   Choosing any two edges from each ~$F_i$, we obtain an~$n$-vertex graph with at least~$(1+\delta)n$ edges, where $\delta=2\alpha-1>0$. Then~\refT{thm:girthsparsegraph} implies that there is a cycle of length at most $C\log n$ for some positive $C(\alpha)$.
   If such a cycle is not rainbow, we can replace two edges in the same edge set~$F_i$ by the other edge in the triangle~$F_i$ to get a shorter cycle. Do it repeatedly until we obtain a rainbow cycle, which is of length at most~$C\log n$. 
\end{proof}

The  next example, the crown-like graph, shows that 
the condition $\alpha > \frac{1}{2}$ is necessary, namely for $\alpha=\frac{1}{2}$ the rainbow girth can be linear in~$n$, not logarithmic. 

\begin{example}
    Let $m=\lfloor\frac{1}{2}n \rfloor$. Let~$K$ be a cycle on  $m$ vertices with edges $e_1,\ldots ,e_m$. Let $v_1, \ldots ,v_m$ be distinct vertices not on~$K$, and let~$F_i$ be the triangle with vertex set $e_i\cup\{v_i\}$. 
    The rainbow girth is~$m$.
    \end{example}

The following theorem implies that the~$\log n$ bound in~\refT{thm:triangle} is the right order of magnitude. The following is a fine-tuned version of Theorem \ref{Omegalogn} from the introduction:

\begin{theorem}\label{thm:nosmallcycle}
    For any~$\alpha>0$, there exists a constant $c>0$ such that for any integer~$n$, there exists an $n$-vertex graph~$G$ formed by at least $\alpha n$ edge-disjoint triangles such that any rainbow cycle in~$G$ has length at least~$c\log n$. 
\end{theorem}

We use two probabilistic tools, the inequalities of Chernoff  and Markov.
\begin{theorem}[Chernoff]\label{thm:Chernoff}
Let $X$ be a binomial random variable $\Bin(n,p)$. For any $0<\epsilon<1$, we have
\[\cPr(X\ge (1+\epsilon)\E X)\le \exp(-\epsilon^2\E X/3). \]
%For any $t\ge 0$, we have
%\[\cPr(X\ge \E X + t)\le \exp\Big(-\frac{t^2}{2(\E X+t/3)}\Big). \]
\end{theorem}

\begin{theorem}[Markov]\label{thm:Markov}
Let $X$ be a non-negative random variable. %For any $0<\epsilon<1$, we have
%\[\cPr(X\ge (1+\epsilon)\E X)\le \exp(-\epsilon^2\E X/3). \]
For any $t>0$, we have
\[\cPr(X\ge t)\le \E X/t. \]
\end{theorem}

\begin{proof}[Proof of~\refT{thm:nosmallcycle}]

Let $p:=\frac{25\alpha}{n^2}$.
    Denote by $G^{(3)}(n,p)=:H$  the system of triples in which each element  of~$\binom{[n]}{3}$ is included independently with probability~$p$.   The example proving the theorem will be the set of triangles induced by the triples in~$G^{(3)}(n,p)$, with some triples removed.

Here are the details.  We have
\[\E |H|=\binom{n}{3}p\ge 4\alpha n.\] 
Chernoff's inequality yields
\begin{equation*}
    \cPr(|H|\le 3\alpha n)\le \cPr(|H|\le 0.9\cdot \E |H|)=o(1).
\end{equation*}

Let~$\cA$ be the event $\{H: |H|\ge 3\alpha n\}$.
Then 
\begin{equation}\label{eq:numberoftriples}
    \cPr(\cA)=1-o(1)
\end{equation}

%with high probability\footnote{An event holds with high probability if the probability of that event tends to 1 as~$n$ tends to infinity.} the number of triples is at least $0.9\binom{n}{3}p\ge 2.5\alpha n$.

Let
\[Y:=|\{(A_1,A_2): A_i\in H \text{ for $i=1,2$} ,A_1\neq A_2\text{ and }|A_1\cap A_2|=2\}|\]
be the number of pairs of distinct triples in~$H$ that intersect at two vertices.

Then  
\[\E Y\le \binom{n}{3}\cdot 3\cdot n p^2 =o(n), \]
as there are at most~$\binom{n}{3}$ ways to choose~$A_1\in\binom{[n]}{3}$,  3 ways to choose a pair $\Pi$ of vertices in the intersection, and then at most~$n$ ways to complete $\Pi$ to the triple ~$A_2\in\binom{[n]}{3}$, and the probability that both~$A_1,A_2$ are in~$H$ is~$p^2$. Then
by Markov's inequality, we have
\[\cPr(Y\ge \alpha n)=o(1). \]
Let $\cB$ be the event that $Y\le \alpha n$. By the above 
\begin{equation}\label{eq:pairsoftriples}
    \cPr(\cB)=1-o(1).
\end{equation}

%By deleting at most one triple of such pair from~$H$, we get~$H^*$ that contains at least~$2\alpha n$ triples such that any two of them intersect in at most 1 vertex.

Given a 3-graph~$F$ on~$[n]$, let $E^{(2)}(F):=\{e\in \binom{[n]}{2}: e\subseteq  A \text{ for some $A\in F$}\}$.
Note that for a fixed~$e\in\binom{[n]}{2}$,
\[\cPr(e\in E^{(2)}(H))=1-\cPr(e\not\subseteq A \text{ for any $A\in H$})=1-(1-p)^{n-2}=:q. \]
By Bernoulli's inequality, we have~$(1-p)^{n-2}\ge 1-(n-2)p$. Therefore
\begin{equation}\label{eq:boundonq}
    q\le (n-2)p\le \frac{25\alpha}{n}.
\end{equation}
Let $C$ be a $k$-cycle in~$K_n$ with edges~$e_1,\dots,e_k$. 
We say that $C$ is {\em distinguishable} if  each $e_i\subseteq  A_i$ for some $A_i\in H$ and $e_j\not\subseteq  A_i$ if $i\neq j$.

We consider the probability that~$C\subseteq E^{(2)}(H)$ and~$C$ is distinguishable. 
We shall bound this probability from above by~$q^k$.
Indeed, for $e_1,\dots, e_k$ we define the event $\cS_{e_i}=\cS_{e_i}(e_1,\dots,e_{i-1})$ as
\begin{align*}
   \cS_{e_i}:=\{\text{there exists } A_i\in\binom{[n]}{3} \text{ with }e_i\subseteq A_i \text{ and $e_j\not\subseteq A_i$ for $j<i$ such that } A_i\in H\}.
\end{align*}

Then we have 
\begin{align*}
    &\cPr(C\subseteq E^{(2)}(H) \text{ is distinguishable})\\
 \le &\cPr(\cap_{i=1}^k \cS_{e_i})=  \prod_{i=1}^{k}\cPr(\cS_{e_i}\mid \cap_{j=1}^{i-1}\cS_{e_j})
\le  (1-(1-p)^{n-2})^k=q^k, 
\end{align*}
where the second to last inequality is because there are at most~$n-2$ many~$A\in\binom{[n]}{3}$ satisfying that~$e_i\subseteq A$ and $e_j\not\subseteq A$ for all $j<i$, thus $\cPr(\neg\cS_{e_i}\mid \cap_{j=1}^{i-1}\cS_{e_j})\ge (1-p)^{n-2}$.

Let~$X_k$ be the number of distinguishable cycles of length~$k$ in~$E^{(2)}(H)$. With a look at~\eqref{eq:boundonq}, we have
\begin{equation*}
    \E X_k\le n^k q^k\le (25\alpha )^k \le n^{1/2}
\end{equation*}
for $k\le c(\alpha)\log n$ and $c>0$ small enough, therefore
\begin{equation*}%\label{eq:EsumXk}
    \sum_{k=3}^{\lfloor c\log n \rfloor} \E X_k= o(n).
\end{equation*}
Then by Markov's inequality, we have 
\[\cPr( \sum_{k=3}^{\lfloor c\log n \rfloor} X_k\ge \alpha n )=o(1)\]
so that
\begin{equation}\label{eq:distinguishablecycles}
    \cPr(\cC)=1-o(1)
\end{equation}
for the event
$\cC:=\{ \sum_{k=3}^{\lfloor c\log n \rfloor} X_k\le \alpha n\}$.

Combining~\eqref{eq:numberoftriples},~\eqref{eq:pairsoftriples}, and~\eqref{eq:distinguishablecycles}, we take~$H$ when~$\cA\cap\cB\cap\cC$ holds, which holds with probability $1-o(1)$. From~$H$, we remove at most one triple in the pairs counted by~$Y$ to get~$H_1$ so that the triples in~$H_1$ intersect with each other in at most one vertex. In particula, each~$e\in E^{(2)}(H_1)$ is contained in exactly one~$A\in H_1$. Then~$\cA\cap \cB$ implies that 
\[|H_1|\ge 3\alpha n-\alpha n \ge 2\alpha n.\] 
View each triple in~$H_1$ as a triangle in~$E^{(2)}(H_1)$. The above observation confirms that the triangles are edge-disjoint. If there is a rainbow cycle in~$E^{(2)}(H_1)$ with edges~$e_1,\dots, e_k$, then $e_i\subseteq A_i\in H_1$, the rainbow property and the fact that there is exactly one triple in~$H_1$ contains an edge in~$E^{(2)}(H_1)$ implies that the cycle is distinguishable. For each rainbow cycle of length at most~$c\log n$ in~$H_1$, in order to destroy the rainbow cycle, we choose at most one edge~$e$ and remove the triple~$A\supseteq e$ from~$H_1$ to get~$H_2$. As~$E^{(2)}(H_1)\subseteq E^{(2)}(H)$, the event~$\cC$ implies that we only need to remove at most~$\alpha n$ triples. Therefore
\[ |H_2|\ge |H_1|-\alpha n\ge \alpha n.  \]
Let~$G:=E^{(2)}(H_2)$. Then~$G$ is a graph formed by at least~$\alpha n$ edge-disjoint triangles without rainbow cycles of length less than~$c\log n$. This completes the proof. 
\end{proof}

\section{Rainbow girth when $\max |F_i| =2$}\label{sec:max2}

For $\cf=(F_1,\dots, F_m)$ a family of subsets of~$E(K_n)$, recall that
\[\psi(\cf)=\sum_{1\le i \le m}\frac{1}{|F_i|}.\]

%\begin{conjecture}\label{rainboeconj}
%Any family  $\cf=(F_1,\dots,F_n)$ of $n$ subsets of  $E(K_n)$ has a rainbow cycle  of size at most $2 %\psi(\cf)$.
%\end{conjecture}

\begin{theorem}\label{mainrainbow}
Let $\cf=(F_1,\dots,F_n)$ be a family of subsets of $E(K_n)$ such that $1\le |F_i| \le 2$. Then $rgirth(\cf) \le \lceil\psi(\cf)\rceil$.

\end{theorem}

In~\cite{DDFGGHMM2021} \refCon{con:generalCHC} was proved when $|F_i|=2$ for all $i$. Theorem \ref{mainrainbow} is a generalization to the case in which some of the sets $F_i$ are singleton sets.

The theorem is easily seen to be equivalent to:

\begin{theorem}\label{mainrainbow2}
Let $\cf=(F_1,\dots,F_n)$ be a family of subsets of $E(K_n)$ such that $1\le |F_i| \le 2$. Assume $p$ sets are of size $1$, and $n-p$ are of size $2$. Then $rgirth(\cf) \le \lceil\frac{n+p}{2}\rceil$.
\end{theorem}

We will refer to the  edges in $F_i$ as colored by  color $i$.

\begin{proof}
We may assume that the sets $F_i$ are disjoint, or else there is a rainbow digon (cycle of length $2$). 
The case where all the sets $F_i$ are of size 2 was proved in \cite{DDFGGHMM2021}. Thus we may assume  $|F_1|=1$. Let $F_1=\{e\}$.

%Assume, first, that there is no vertex connected in one color to both endpoints of $e$. Then after contracting $e$ to a single vertex $h$ we have $p-1$ color classes of size 1 and by the induction hypothesis there exists a rainbow cycle of length $\lceil \frac{n-1+p-1}{2}\rceil$. Opening up this cycle (which may require using $e$ in the new cycle) yields a cycle of length at most $\lceil \frac{n+p}{2}\rceil$, proving the desired result. 

%So, we may assume that there is a vertex $x_1$ connected by edges in (say) $F_2$ to both endpoints of $e$. Let $H_1$ be the triangle on $V(e) \cup \{x_1\}$. 

%This is the first step in a recursive procedure, 

We construct a subgraph~$H$ of~$G$ recursively as follows. Let $H_0=\{e\}$. At each step $i$,  $H_i$ is constructed by adding to $H_{i-1}$ a vertex $x_i \notin V(H_{i-1})$ and two edges $x_ia_i,x_ib_i \notin E(H_{i-1})$ such that $a_i,b_i \in V(H_{i-1})$ and  $x_ia_i,x_ib_i$ are colored by the same color $i$. We stop at step $i=t$ when there are no such two edges to add, and we let $H=H_t$.

For two vertices $u,v\in V(G)$ let $\dist_{r,G}(u,v)$ denote the {\em rainbow distance} of $u,v$, that is, the minimum length (number of edges) of a rainbow path in~$G$ connecting $u,v$.  
For a subgraph~$G'$ of~$G$ let the {\em rainbow diameter} of~$G$ be defined as  $rd(G'):=\max_{u,v\in V(G')}\dist_{r,G'}(u,v)$. We omit the subscript~$G$ in~$dist_r$ and it should be clear according to the context.

\begin{claim}\label{claim:diaHi}
$rd(H_i)\le \frac{i}{2} +1$, and if~$i$ is even, except at most one pair of vertices  $u_i,v_i \in V(H_i)$, for any pair of vertices $u,v\in V(H_i)$, we have $\dist_r(u,v) \le \frac{i}{2}$.
\end{claim}
\begin{proof}
If $i\in \{0,1\}$ the claim is trivial. We proceed by induction on $i$. 

Suppose first that $i+1$ is odd.   By the induction hypothesis, there exists at most one pair of vertices $u_i,v_i\in V(H_i)$ such that  $\dist_r(u_i,v_i) = \frac{i}{2} +1$ and for any other pair of vertices $u,v\in V(H_i)$, $\dist_r(u,v) \le \frac{i}{2}$.
We have to show that for every $y,z \in V(H_{i+1})$, $\dist_r(y,z) \le \lfloor \frac{i+1}{2} +1 \rfloor = \frac{i}{2} +1$. If $y,z\in V(H_{i})$ we are done. Suppose $z=x_{i+1}$. 
If $y\notin \{u_i,v_i\}$ there is a rainbow path from $a_{i+1}$ to $y$ of length at most $\frac{i}{2}$ and thus there is a rainbow path from $x_{i+1}$ to $y$ of length at most $\frac{i}{2}+1$. 
If $y \in \{u_i,v_i\}$, say $y=u_i$, then either $a_{i+1}\neq v_i$ or $b_{i+1}\neq v_i$. In both cases there exists a rainbow path from $x_{i+1}$ to $y$, through $a_{i+1}$ or $b_{i+1}$ respectively, of length at most $\frac{i}{2}+1$. 

Assume now that~$i+2$ is even. 
By the induction hypothesis, there exists at most one pair $u_i,v_i\in V(H_i)$ such that $\dist_r(u_i,v_i) = \frac{i}{2} +1$ and  any other pair of vertices in $V(H_i)$ is of rainbow  distance at most $\frac{i}{2}$.
We have to show that there is at most one pair $u_{i+2},v_{i+2}\in V(H_{i+2})$ such that  $\dist_r(u_{i+2},v_{i+2}) = \frac{i}{2} +2$ and  any other pair of vertices in $V(H_{i+2})$   is of rainbow distance at most  $\frac{i}{2}+1$.

We split into two cases. 
\medskip

{\bf Case 1.} $x_{i+1} \notin \{a_{i+2},b_{i+2}\}$.  

Choose $u_{i+2}=x_{i+1}, v_{i+2}=x_{i+2}$. We claim that $\dist_r(u_{i+2},v_{i+2}) \le \frac{i}{2} +2$. Indeed, by the induction hypothesis we can choose a vertex  $u\in \{a_{i+1}, b_{i+1}\}$ and $v\in \{a_{i+2}, b_{i+2}\}$ such that $\dist_r(u,v) \le \frac{i}{2}$ by the fact that $a_{i+1},b_{i+1},a_{i+2},b_{i+2} \in V(H_{i})$ and $\{u,v\}\neq \{u_i,v_i\}$, and then adding the edges $x_{i+1}u, x_{i+2}v$ we get a rainbow path between $x_{i+1},x_{i+2}$ of length at most $\frac{i}{2} +2$.
Let $u,v\in V(H_{i+2})$ such that $\{u,v\}\neq \{x_{i+1},x_{i+2}\}$. Since $a_{i+1},b_{i+1},a_{i+2},b_{i+2} \in V(H_{i})$, we have $\dist_r(u,v) \le  \frac{i}{2} +1$ like in the odd case.

\medskip

{\bf Case 2.} $x_{i+1} = a_{i+2}$.  

In this case, either $b_{i+2} \neq u_i$ or $b_{i+2} \neq v_i$. Assume WLOG $b_{i+2} \neq v_i$. Choose $u_{i+2}=x_{i+2}$, $v_{i+2}=v_{i}$. 
Then $\dist_r(u_{i+2},v_{i+2}) \le \lfloor \frac{i+1}{2} +1\rfloor +1=  \frac{i}{2} +2$, and like before, $\dist_r(u,v) \le  \frac{i}{2} +1$ for any other pair of vertices $u,v$. 
\end{proof}

Returning to the proof of the theorem, we proceed by induction on $n$. Contract $H$ into a single vertex $h$ to obtain a new graph $G'$ ($G'$ may have loops). Note that $n':=|V(G')|=n-t-1 $, the number of colors is $n-t-1=n'$ and the number of colors of size 1 is $p'=p-1$. By induction there exists a rainbow cycle $C$ in $G'$ of size at most $\lceil\frac{n'+p'}{2}\rceil = \lceil\frac{n-t+p-2}{2}\rceil= \lceil\frac{n-t+p}{2}\rceil -1$. 

If~$C$ does not use the vertex~$h$, we are done. Otherwise, uncontracting~$h$,~$C$ either remains a cycle (possibly containing a vertex in~$h$) - in this case we are done; or it may become  a path in~$G$, with end vertices $u\neq v\in V(H)$.     
By~\refCl{claim:diaHi}, there is a rainbow path~$P$ in~$H$ connecting~$u$ and~$v$ of size at most $\frac{t}{2}+1$. Note that~$P$ uses colors not appearing in~$C$. Thus~$P+C$ is a rainbow cycle in~$G$ of size 
at most $\lfloor \lceil\frac{n-t+p}{2}\rceil -1 + \frac{t}{2}+1 \rfloor = \lceil\frac{n+p}{2}\rceil$.
This completes the proof of the theorem.
\end{proof}

\begin{corollary}\label{cor:directedp}
    Let $D$ be an $n$-vertex sink-less digraph. Assume $p$ vertices have out-degree $1$. Then 
$girth(D)  \le \lceil\frac{n+p}{2}\rceil$.
\end{corollary}

\begin{remark}
In \cite{shen1} (Theorem 1) a slightly weaker result was proved: $dgirth(D) \le \lceil\frac{n+p+1}{2}\rceil$.
\end{remark}

\begin{proof}[Proof of~\refC{cor:directedp}]
    For each vertex~$v$ of~$D$ that has out-degree more than 2, we remove some arbitrary edges to make~$v$ have out-degree exactly 2. Then in the resulting digraph~$D'$, there are~$p$ vertices of out-degree 1 and~$n-p$ vertices of out-degree 2, and we have $girth(D)\le girth(D')$. Therefore it is enough to prove that $girth(D') \le \lceil\frac{n+p}{2}\rceil$. 
    While by the construction to explain why~\refCon{con:generalCHC} generalizes CHC in~\refS{sec:triangles}, we reduce the problem into rainbow undireceted version with~$p$ stars of size 1 and $n-p$ stars of size 2. 
    %We construct an $n$-vertex graph~$G$ as follows: we set $V(G)=V(D')$ and let the outward edges of~$D'$ at each vertex~$v$ be a set of edges~$F_v$, forgetting the directions. Then we get~$n$ sets of edges, $p$ of which are of size~$1$ and $n-p$ of which are of size 2. (Actually each $F_v$ is a star.)  Note that a rainbow cycle in~$G$ correspond to a directed cycle in~$D'$, therefore 
    Therefore we complete the proof by applying~\refT{mainrainbow2}. 
\end{proof}

\begin{corollary}
For a family $\cf=(F_1,\dots, F_n)$ of subsets of~$E(K_n)$ satisfying $1\le |F_i|\le 2$ for all $1\le i\le n$, we have
$rgirth(\cf) \le 2\sum_{1\le i \le n}\frac{1}{|F_i|+1}$.
\end{corollary}
\begin{proof}
    Applying~\refT{mainrainbow2}, we have $rgirth(\cf) \le \lceil\psi(\cf)\rceil = \lceil\frac{n+p}{2}\rceil$, where~$p$ is the number of sets in~$\cf$  of size 1. Note that $\lceil\frac{n+p}{2}\rceil\le \frac{n+p}{2}+\frac{1}{2}$, which is at most $2(\frac{p}{2}+\frac{n-p}{3})=2\sum_{1\le i \le n}\frac{1}{|F_i|+1}$ when $p\le n-3$. Furthermore, for $p=n-2$ or~$n$, $\lceil\psi(\cf)\rceil = \lceil\frac{n+p}{2}\rceil=\psi(\cf)\le 2\sum_{1\le i \le n}\frac{1}{|F_i|+1}$. And in the remaining case $p=n-1$, we have $rgirth(\cf)\le n-1 \le \psi(\cf)\le 2\sum_{1\le i \le n}\frac{1}{|F_i|+1}$.
\end{proof}

\end{document}